\documentclass[11pt,oneside]{amsart}
\usepackage{amsopn}
\usepackage{amsmath,amsthm,amssymb}
\usepackage{prettyref}

\newcommand{\nc}{\newcommand}

 \nc{\aff}{\mathfrak{aff} } \nc{\bb}{\mathfrak{b} }
\nc{\cc}{\mathfrak{c} }  \nc{\dd}{\mathfrak{d} }
 \nc{\ggo}{\mathfrak{g} }
 \nc{\hh}{\mathfrak{h} }  \nc{\ii}{\mathfrak{i} }
 \nc{\jj}{\mathfrak{j} }  \nc{\kk}{\mathfrak{k} }
\nc{\mm}{\mathfrak{m} }   \nc{\nn}{\mathfrak{n} }
\nc{\pp}{\mathfrak{p} }  \nc{\rr}{\mathfrak{r} } \nc{\sg}{\mathfrak{s} }
 \nc{\sog}{\mathfrak{so} }  \nc{\spg}{\mathfrak{sp} }
 \nc{\sug}{\mathfrak{su} }  \nc{\slg}{\mathfrak{sl} }
 \nc{\tg}{\mathfrak{t} }  \nc{\uu}{\mathfrak{u} }
 \nc{\vv}{\mathfrak{v} } \nc{\ww}{\mathfrak{w} }
 \nc{\zz}{\mathfrak{z} }

 \nc{\ggob}{\overline{\mathfrak{g}}}

\nc{\glg}{\mathfrak{gl} }

\nc{\pca}{\mathcal{P}} \nc{\nca}{\mathcal{N}}

 \nc{\vp}{\varphi} \nc{\ddt}{\frac{{\rm d}}{{\rm d}t}}
 \nc{\la}{\langle} \nc{\ra}{\rangle}

 \nc{\SO}{{\sf SO}} \nc{\Spe}{{\sf Sp}} \nc{\Sl}{{\sf Sl}}
 \nc{\SU}{{\sf SU}} \nc{\Or}{{\sf O}} \nc{\U}{{\sf U}}
 \nc{\Gl}{{\sf Gl}} \nc{\Se}{{\sf S}} \nc{\Cl}{{\sf Cl}}
 \nc{\Spin}{{\sf Spin}} \nc{\Pin}{{\sf Pin}}

 \nc{\RR}{{\mathbb R}} \nc{\HH}{{\mathbb H}} \nc{\CC}{{\mathbb C}}
 \nc{\ZZ}{{\mathbb Z}} \nc{\FF}{{\mathbb F}} \nc{\NN}{{\mathbb N}}
 \nc{\GG}{{\mathbb G}} \nc{\JJ}{{\mathbb J}} \nc{\II}{{\mathbb I}}
 \nc{\KK}{{\mathbb K}} \nc{\DD}{{\mathbb D}}

 \nc{\ad}{\operatorname{ad}} \nc{\Ad}{\operatorname{Ad}}
 \nc{\coad}{\operatorname{coad}} \nc{\ct}{\operatorname{T}}
 \nc{\rank}{\operatorname{rank}} \nc{\Irr}{\operatorname{Irr}}
 \nc{\End}{\operatorname{End}} \nc{\Aut}{\operatorname{Aut}}
 \nc{\Inn}{\operatorname{Inn}} \nc{\Der}{\operatorname{Der}}
 \nc{\Dera}{\operatorname{Dera}} \nc{\Auto}{\operatorname{Auto}}
 \nc{\GL}{\operatorname{GL}}
 \nc{\SL}{\operatorname{SL}}

 \theoremstyle{plain}
 \newtheorem{Theorem}{Theorem}[section]
 \newtheorem{Proposition}[Theorem]{Proposition}
 \newtheorem{Corollary}[Theorem]{Corollary}
 \newtheorem{Lemma}[Theorem]{Lemma}

 \newenvironment{Proof}{\noindent \emph{Proof.}}{\hfill $\blacksquare$}

 \newtheorem*{Remark}{Remark}
 \newtheorem{Definition}[Theorem]{Definition}
 
 \newtheorem{Example}[Theorem]{Example}

 \newcommand{\R}{\mathbb R}
\newcommand{\N}{\mathbb N}
\newcommand{\Z}{\mathbb Z}
\newcommand{\C}{\mathbb C}

\newcommand{\mg}{\mathfrak g }

\newcommand{\mn}{\mathfrak n }
\newcommand{\mz}{\mathfrak z }
\newcommand{\mv}{\mathfrak v }
\newcommand{\mh}{\mathfrak h }
\newcommand{\ma}{\mathfrak a }
\newcommand{\mgg}{\mathfrak g }

\newcommand{\tmn}{\tilde{\mathfrak n}}

\newcommand{\wtV}{\widetilde V }

\newcommand{\tomega}{ \tilde{\omega} }

\newcommand{\bil}{\left\langle \;,\;\right\rangle}
\newcommand{\lela}{\left \langle}
\newcommand{\rira}{\right \rangle}
\newcommand{\lra}{\longrightarrow}

\begin{document}

\title{Symplectic structures on nilmanifolds: an obstruction for its existence}

\author{Viviana J. del Barco
}
\address{V. del Barco: ECEN-FCEIA, Universidad Nacional de Rosario, Pellegrini 250, 2000 Rosario, Argentina.
}

\email{delbarc@fceia.unr.edu.ar}

\thanks{Supported by a fellowship from CONICET and research grants from SCyT-UNR, Secyt-UNC}





\begin{abstract} 
In this work we introduce an obstruction for the existence of symplectic structures on nilpotent Lie algebras. Indeed, a necessary condition is presented in terms of the cohomology of the Lie algebra. Using this obstruction we obtain both positive and negative results on the existence of symplectic structures on a large family of nilpotent Lie algebras. Namely the family of nilradicals of minimal parabolic subalgebras associated to the real split Lie algebra of classical complex simple Lie algebras.
\medskip



\end{abstract}

\maketitle

\section{Introduction}

A nilmanifold  is an homogeneous manifold $M=\Gamma\backslash N$ where $N$ is a simply connected nilpotent Lie group  and $\Gamma$ is co-compact discrete subgroup of $N$. For these compact manifolds the natural map from $H_{dR}^i (\mn)$, $\mn$ the Lie algebra of $N$, to the de Rham cohomology group $H^i(M,\R)$ is an isomorphism for all $0 \leq i \leq 2n$, as showed by Nomizu in \cite{NO}.

In particular this implies that any symplectic structure on a nilmanifold is cohomologous to an invariant one. Thus to solve the problem of existence of symplectic structures on the nilmanifold $\Gamma \backslash N$ reduces to find a non-degenerate closed 2-form $\omega$ on the Lie algebra $\mn$; if it exists $\mn$ is called a symplectic Lie algebra. Here we work from this Lie algebra point of view.

The goal of this work is to prove that every symplectic nilpotent Lie algebra has a certain non-zero component on its cohomology. Actually, the intermediate cohomology of a Lie algebra $\mn$ (concept presented by the author in \cite{dB2}) is used in Theorem \ref{teo:teo3} to give a necessary condition for $\mn$ to admit a symplectic structure. As an application, we study the validity of this property on a particular subfamily of nilpotent Lie algebras. 

\smallskip

Benson and Gordon in \cite{BE-GO} proved that the Hard Lefschetz Theorem fails for any symplectic non-abelian nilpotent Lie algebra. In order to show this, they deduce some general structure results of symplectic nilpotent Lie algebras. Nevertheless there are not many general conditions to determine whether a given nilpotent Lie algebra is symplectic or not.

Until now there are known all the symplectic nilpotent Lie algebras up to dimension $6$ (see \cite{BO-GO,SA1} for instance) and this list is mostly build-up by studying case by case. But the lack of a full classification of real nilpotent Lie algebras of dimension $\geq 8$ makes this method non-feasible to find the symplectic Lie algebras in greater dimensions.

Moreover, several authors studied the problem on different subfamilies of nilpotent Lie algebras. For example, the classification of symplectic filiform Lie algebras, which are Lie algebras $\mn$ of nilpotency index $k=\dim \mn-1$, is given in \cite{MI}. Moreover, in \cite{DO-TI} the authors work with Heisenberg type nilpotent Lie algebras. Among nilpotent Lie algebras associated with graphs, a complete description of the symplectic ones can be made in terms of the corresponding graph \cite{PO-TI}. The full classification of the symplectic free nilpotent Lie algebras is done in \cite{dB}. 

In this context, the aim of this work is to contribute with a better understanding of the structure of symplectic nilpotent Lie algebras. Its organization is as follows. Section 2. is devoted to an introduction to the intermediate cohomology of nilpotent Lie algebras and the development of the properties that will be used later on the presentation. In Section 3. we study the relationship between symplectic structures and intermediate cohomology. This leads us to a necessary condition for a nilpotent Lie algebra to admit a symplectic structure. We notice that this condition is not sufficient in general.

In Section 4 we restrict ourselves to the study of the existence of symplectic structures in the family of nilradicals of minimal parabolic subalgebras associated to the real split Lie algebras corresponding to complex classical simple Lie algebras. We prove that the obstruction in Theorem \ref{teo:teo3} is also sufficient for that family. 
This allows us to obtain both positive and negative results about the existence of symplectic structures in this case.

Recall that for the nilpotent complex case Kostant in \cite{KO} describes the Lie algebra cohomology groups of the nilradicals of Borel subalgebras for any irreducible representation as a direct sum of one dimensional modules of multiplicity one. The real version of his description was recently given in \cite{SI}. Here we also use a decomposition of the cohomology groups but the summands are not in one to one correspondence with those of neither Kostant (in the complex version) nor {\v{S}}ilhan in the real case.

\section{Intermediate Cohomology of nilpotent Lie algebras}
The concept of intermediate cohomology of nilpotent Lie algebras and a deep study of its properties were analyzed by the author in \cite{dB2}. For completeness of this work we give here a brief introduction to this cohomology by quickly reviewing its definition and the properties that will be used later.\smallskip

Let $\mgg$ denote a real Lie algebra. The central descending series of $\mgg$, $\{\mgg^i\}$ for all $i\geq 0$,  is given by
$$\mgg^0=\mgg,\quad  \mgg^i=[\mgg,\mgg^{i-1}],\;\; i\geq 1.
$$
A Lie algebra $\mg$ is $k$-step nilpotent if $\mgg^k=0$ and $\mgg^{k-1}\neq 0$; this number $k$ is called the nilpotency index of $\mg$. Nilpotent Lie algebras will be denoted by $\mn$. Abelian Lie algebras are 1-step nilpotent. Moreover, $2$-step nilpotent Lie algebras verify $\mn^1\subseteq \mz(\mn)$, where $\mz(\mn)$ denotes the center of $\mn$.

\smallskip

The Chevalley-Eilenberg complex of a Lie algebra $\mgg$ of dimension $m$ is
\begin{equation}\label{eq:complex}0\lra\R \lra \mgg^*\stackrel{d_1}{\lra}\Lambda^2\mgg^*\stackrel{d_2}{\lra}\ldots \ldots\stackrel{d_{m-1}}{\lra} \Lambda^m\mg^*\lra 0\;\,.\end{equation}
We identify the exterior product $\Lambda^p\mgg^*$ with the space of skew-symmetric $p$-linear forms on $\mgg$, thus each differential $d_p:\Lambda^{p}\mgg^*\lra \Lambda^{p+1}\mgg^*$ is defined by:
$$ d_p c\,(x_1,\ldots,x_{p+1})=\sum_{1\leq i<j\leq p+1}(-1)^{i+j-1}c([x_i,x_j],x_1,\ldots,\hat{x_i},\ldots,\hat{x_j},\ldots,x_{p+1}). $$

The first differential $d_1$ coincides with the dual mapping of the Lie bracket $[\,,\,]:\Lambda^2\mgg\lra \mgg$ and the collection of $d_p$ is a derivation of the exterior algebra $\Lambda^*(\mgg^*)$. We will denote $d$ instead of $d_p$ independently of $p$.  

The cohomology of $(\Lambda^*\mgg^*,d)$ is called the Lie algebra cohomology of $\mgg$ (with real coefficients) and it is denoted by $H^*(\mgg,\R)$ and more often as $H^*(\mgg)$ if there is no place to confusion. For nilpotent Lie algebras $H^1(\mn)\cong \mn/\mn^1$ and $\dim H^2(\mn) \geq 2$ \cite{DI2}.
\medskip

When the Lie algebra is nilpotent, a filtration of the cochain complex in Eq. (\ref{eq:complex}) arises in the following manner.
Consider the subspaces of $\mn^*$ defined by Salamon in \cite{SA1}
\begin{equation} \label{eq:eq53} 
V_0=0 \qquad \qquad V_i=\{ \alpha \in \mn^*: d\alpha \in \Lambda^2V_{i-1}\}\qquad i\geq 1.
\end{equation}
Then $V_0 \subseteq V_1\subseteq \cdots \subseteq V_i\subseteq \cdots\subseteq \mn^* $ and $V_i$ is the annihilator of $\mn^i$, the $i$th-ideal in the central descending series; that is $V_i=(\mn^i)^\circ$. In particular, $\mn$ is a $k$-step nilpotent Lie algebra if and only if $V_k=\mn^*$ and $V_{k-1}\neq \mn^*$.

\medskip
Suppose $\mn$ is a $k$-step nilpotent Lie algebra of dimension $m$, then for any $q=0,\ldots,m$, the space of skew symmetric $q$-forms $\Lambda^q\mn^*$ is filtered since
\begin{equation}\label{l4}
0=\Lambda^q V_0 \subsetneq \Lambda^q V_1 \subsetneq \ldots \subsetneq \Lambda^q V_{k-1} \subsetneq \Lambda^q V_k=\Lambda^q \mn^*.
\end{equation}

In addition each of these subspaces is invariant under the differential, therefore
 \begin{equation}\label{eq7}F^pC^*: 0 \lra \R \lra V_{k-p} \lra \Lambda^2 V_{k-p} \lra \cdots \lra \Lambda^m V_{k-p} \lra 0 \end{equation}
  is a subcomplex of the Chevalley-Eilenberg complex for each fixed $p$ and $\{F^pC^*\}_{p\geq 0}$ constitutes a filtration of the complex in Eq. \prettyref{eq:complex}.

As any filtration of a cochain complex, $\{F^pC^*\}_{p\geq 0}$ gives rise to a spectral sequence $\{E_r ^{p,q}(\mn)\}_{r\geq 0} ^{p,q \in \Z}$. In this case, this spectral sequence always
converges to the Lie algebra cohomology of $\mn$ (see \cite{dB2} and references therein). 
In particular this implies that each cohomology group $H^i(\mn)$ can be written as a direct sum of the limit terms of the spectral sequence. Namely
\begin{equation}
\label{eq:eq8}H^i(\mn)\cong \bigoplus _{p+q=i} E_\infty^{p,q}(\mn) \quad \text{ for all } i=0,\ldots, m.
\end{equation} 

This way of describing the cohomology groups as a sum of smaller spaces suggests us the following definition.

\begin{Definition} Let $\mn$ be a nilpotent Lie algebra of dimension $m$. Then, for each $i=0,\ldots, m$, the intermediate cohomology groups of degree  $i$ of $\mn$ are the vector spaces $E^{p,q}_\infty(\mn)$ with $p+q=i$.
\end{Definition}

Notice that for each $i=0,\ldots, m$ there is a finite amount of non-zero intermediate cohomology groups of degree $i$.

Each intermediate cohomology group can be described using the Lie algebra differential restricted to the subspaces in the filtration: \begin{equation}\label{eq6}E_\infty^{p,q}(\mn)\cong\frac{\{x \in \Lambda^{p+q}V_{k-p}:dx =0\}}{d(\{x\in\Lambda^{p+q-1}\mn^*:dx\in \Lambda^{p+q}V_{k-p}\})+\{x\in\Lambda^{p+q}V_{k-p-1}:dx=0\}}. \end{equation}

If a nilpotent Lie algebra $\mn$ can be decomposed as a direct sum of a one dimensional ideal $\R$ and a nilpotent Lie algebra of dimension one less than $\mn$, a similar formula to the K\"unneth formula can be stated for the intermediate cohomology.

\begin{Theorem}[\cite{dB2}]\label{teo:teo24} Let $\mn$ be a $k$-step nilpotent Lie algebra which can be decomposed as a direct sum of ideals $\mn=\R\oplus\mh$. Then $\mh$ is $k$-step nilpotent and for all $0\leq r\leq \infty$ it holds
\begin{enumerate}
\item $E_r^{p,-p}(\mn)=0$ for all $p=0,\ldots,k-2$ and $E_r^{k-1,1-k}(\mn)\cong\R$.
\item $E_r^{k-1,2-k}(\mn)\cong E_r^{k-1,2-k}(\mh)\oplus \R$,
\item $E_r^{p,1-p}(\mn)\cong E_r^{p,1-p}(\mh)$ \hspace{0.5cm} if $p \leq k-2$,
\item $E_r^{p,q}(\mn)\cong E_r^{p,q}(\mh)\oplus E_r^{p,q-1}(\mh)$  \hspace{0.5cm} if $p+q\geq 2$.
\end{enumerate}
\end{Theorem}

Throughout an inductive procedure the next result follows.

\begin{Corollary} \label{cor:cor1}Suppose $\mn$ is a non-abelian nilpotent Lie algebra. Then $E_\infty^{0,2}(\R^s\oplus \mn)= E_\infty^{0,2}(\mn)$ for any $s\geq 0$.
\end{Corollary}

\section{Symplectic structures and the $E_\infty^{0,2}$ intermediate cohomology group}

A symplectic structure on a differentiable manifold $M$ is a differentiable closed 2-form $\Omega$ that is non-singular at every point of $M$. Not every manifold admits such a structure. For example, it is well known that compact manifolds having zero second de Rham cohomology group do not admit symplectic structures.  When $M$ is a nilmanifold this criteria is useless since $H^2_{dR}(M)\cong H^2(\mn)$, always non-zero for nilpotent Lie algebras (\cite{DI2}). And yet there exists non-symplectic nilmanifolds. We present here an adapted version of this criteria that can be used to determine non-existence of symplectic structures on nilmanifolds.

Recall that a Lie algebra is symplectic if it admits a skew-symmetric bilinear form $\omega$ which is both closed and is non-degenerate. In Theorem \ref{teo:teo3} we prove that there is a close relationship between the existence of symplectic structures on a nilpotent Lie algebra $\mn$ and its intermediate cohomology group $E_\infty^{0,2}(\mn)$. 
 From this, a new general obstruction for the existence of these structures on nilpotent Lie algebras is deduced. 
\smallskip

Let $\mn$ be a nilpotent Lie algebra and $\omega$ a symplectic structure on $\mn$. Consider an element $z\notin \mn$ and define the bracket $[[\,,\,]]$ in $\tmn=\mn\oplus \R z$ from that one  $[\,,\,]$ in $\mn$ as follows:
\begin{equation}
\left[\left[x,y\right]\right]=[x,y]+\omega(x,y)\,z,\qquad
[[x,z]]=0 \qquad \mbox{ for all }\; x,y\,\in\mn .\nonumber
\end{equation}
 Notice that $\tmn$ is the central extension of $\mn$ by $\omega$; this  extension was considered in \cite{BE-GO} to describe symplectic nilmanifolds as subquotients of coadjoint orbits.
\begin{Lemma}
Let $\mn$ be a nilpotent Lie algebra with a symplectic structure $\omega$ and consider the extension $\tmn$ described above. Then, $\tmn$ is nilpotent and it verifies
\begin{equation}\nonumber
\dim\, E_\infty^{0,2}(\mn)=\dim\,E_\infty^{1,1}(\tmn)+1.
\end{equation}
In particular $E_\infty^{0,2}(\mn)\neq 0$.\end{Lemma}
\begin{Proof}
It is easy to check that $\tmn$ is nilpotent and has its center spanned by $z$. So it is possible to construct on each $\mn^*$ and $\tmn^*$ the filtration described in the previous section.

Denote by $V_i$ y $\wtV_i$  the subspaces of $\mn^*$ and $\tmn^*$ respectively defined in Eq. \prettyref{eq:eq53}. Then
$$ \wtV_i=V_i, \;i=0,\ldots,k,\qquad \wtV_{k+1}=\tmn^*.$$ Let $\tilde{d}$ be the differential of $\tmn$ and let $\zeta\in \tmn^*$ be such that $\zeta(z)=1,\; \zeta(\mn)=0$. It holds  $\tilde{d}\zeta=-\omega$.

Moreover, the restriction of $\tilde{d}$ to $\mn^*$ coincides with the differential $d$ of $\mn$. From (\ref{eq6})
$$
E_\infty^{1,1}(\tmn)
=\frac{\{x\in \Lambda^2\mn^*:dx=0\}}{\R\omega+\{x\in\Lambda^2V_{k-1}:dx=0\}} \quad\mbox{ and }  \quad E_\infty^{0,2}(\mn)=\frac{\{x\in\Lambda^2\mn^*:dx=0\}}{\{x\in\Lambda^2V_{k-1}:dx=0\}}.
$$
Since the symplectic form $\omega$ is non-degenerate, $\omega\notin \Lambda^2V_{k-1}$ and the equation above holds.
\end{Proof}

An immediate consequence is:
\begin{Theorem}\label{teo:teo3} If $\mn$ is a nilpotent Lie algebra which verifies $E_\infty^{0,2}(\mn)=0$, then $\R^s\oplus \mn$ does not admit symplectic structures for all $s\geq 0$.
\end{Theorem}

\begin{Proof} If $E_\infty^{0,2}(\mn)= 0$ then $\mn$ is non-abelian (see \cite[3.1 Examples]{dB2}). Then Corollary \ref{cor:cor1} implies that $E_\infty^{0,2}(\mn)=E_\infty^{0,2}(\R^s\oplus \mn)=0$ for any $s\geq  0$. By the previous Lemma, $\R^s\oplus \mn$ does not admit symplectic structures.
\end{Proof}
\medskip

The reciprocal result to that in Theorem \ref{teo:teo3} is not valid in general as the next example shows.


\begin{Example} Let $\mn_{m,3}$ be the free $3$-step nilpotent Lie algebra on $m$ generators. Recall that $\mn_{m,3}=\mathfrak f_m/(\mathfrak f_m)^{3}$ where $\mathfrak f_m$ is the free Lie algebra on $m$ generators. On the one hand, when $m\geq 3$ the Lie algebra $\mn_{m,3}$ does not admit symplectic structures as proved in \cite{dB}. 

On the other hand, $E_\infty^{0,2}(\mn_{m,3})=0$ for all $m$. Indeed, consider the Hall basis $\mathcal B$ of $\mn_{m,3}$ for a set of generators $\{e_1,\ldots,e_m\}$; the elements in $\mathcal B$ have the form
$$e_i,\quad [e_j,e_k],\quad [[e_r,e_s],e_t], \mbox{ for }\;i=1,\ldots,m,\;1\leq k< j\leq m\mbox{ and }\,1\leq s<r\leq m,\,t\geq s.
$$
 These basis were introduced by Hall in \cite{HA} and they are the usual ones to work with when dealing with free Lie algebras.
The dual basis of $\mathcal B$ consists of 1-forms $\alpha^i,\, \alpha^{jk},\alpha^{rst}$ and its differentials, by Maurer-Cartan formulas, are
$$
\left\{
\begin{array}{l}
d\alpha^i=0,\;i=1,\ldots, m,\\
d\alpha^{ij}=-\alpha^i\wedge \alpha^j, \;1\leq j<i\leq m,\\
d\alpha^{ijk}=-\alpha^{ij}\wedge \alpha^k, \;1\leq j<i\leq m,\, k\geq j.
\end{array}
\right.
$$

The filtration in Eq. \prettyref{eq:eq53} of $\mn_{m,3}^*$ is
$$V_1=\mbox{span}\,\{\alpha^i,\;i=1\ldots,m\},\quad V_2=\mbox{span}\,\{\alpha^i,\,\alpha^{jk}, i=1,\dots,m, 1\leq k< j\leq m\}\;\;\; \mbox{ and } \; \;\; V_3=\mn_{m,3}^*.$$

For each  $1\leq j<i\leq m,\, k\geq j$ the 2-form $\;\alpha^{ijk}\wedge \alpha^k$  defines a non-zero element in $E_\infty^{0,2}(\mn_{m,3})$. In fact
$d(\alpha^{ijk}\wedge \alpha^k)=-\alpha^{ij}\wedge \alpha^k\wedge \alpha^k=0$ so  $\alpha^{ijk}\wedge \alpha^k$ is a closed form in $V_1\wedge (\mn^3)^*$ 
but $\alpha^{ijk}\wedge \alpha^k \notin \Lambda^2V_2$. Therefore $E_\infty^{0,2}(\mn_{m,3})\neq 0$. 
\end{Example}

\begin{Remark} Nilpotent Lie algebras associated with graphs admitting symplectic structures were characterized in terms of their graphs by Poussele and Tirao in \cite{PO-TI}. 
 Meanwhile it is possible to show that all 2-step nilpotent Lie algebras associated with graphs verify $E_\infty^{0,2}\neq 0$. 
Therefore in this family the reciprocal result to that one in Theorem \ref{teo:teo3} neither holds.
\end{Remark}

\subsection{$\operatorname{Aut}(\mn)$ action on $E_\infty^{0,2}(\mn)$}

Once it is known that a certain Lie algebra is symplectic, it is interesting to classify its symplectic forms up to equivalence. In the case of symplectic nilpotent Lie algebras the subspace $E_\infty^{0,2}$ is non-zero. What we study here is how this subspace helps to this classification problem. 


 The automorphism group of a Lie algebra $\mg$ is
$$\Aut(\mg)=\{A\in\,\mbox{GL}(\mg)\,:\,[Ax,Ay]=[x,y]\quad \mbox{for all } x,y\in \mg\}.$$
This group acts on $H^2(\mg)$ in the following way: $A\cdot [\omega]=[(A^{-1})^*\omega]$, for all $A\in\Aut(\mg)$ and $[\omega]\in H^2(\mg)$. Here $(A^{-1})^*$ denotes the automorphism of the exterior algebra $\Lambda^*\mn^*$ induced by $(A^{-1})^*:\mg^*\lra \mg^*$. 


When the Lie algebra is nilpotent, its group of automorphisms $\operatorname{Aut}(\mn)$ acts similarly in $E_\infty^{0,2}(\mn)$. Given a closed 2-form $\omega$ in $\Lambda^2\mn^*$ denote with  $[\omega]^{0,2}$ its class as an element of the quotient space
$$\frac{\{x \in \Lambda^{2}\mn^*:dx =0\}}{\{x\in\Lambda^{2}V_{k-1}:dx=0\}}\cong E_\infty^{0,2}(\mn).$$

Any element in $\Aut(\mn)$ preserves the filtration in (\ref{l4}) of $\mn^*$ and in particular $(A^{-1})^*V_{k-1}=V_{k-1}$. This implies that if $\omega_1,\,\omega_2$ are closed 2-forms on $\Lambda^2\mn^*$ with $[\omega_1]^{0,2}=[\omega_2]^{0,2}$ then $[(A^{-1})^*\omega_1]^{0,2}=[(A^{-1})^*\omega_2]^{0,2}$. Therefore, the following action is well defined: 
\begin{eqnarray}\label{action}
\Aut(\mn)\times E_\infty^{0,2}(\mn)&\lra& E_\infty^{0,2}(\mn)\\
(A,[\omega]^{0,2})&\mapsto& A\cdot[\omega]^{0,2}=[(A^{-1})^*\omega]^{0,2}.\nonumber
\end{eqnarray}

\begin{Proposition} For any nilpotent Lie algebra $\mn$ the map
$$p:H^2(\mn)\lra E_\infty^{0,2}(\mn) ,\quad[\omega]\mapsto [\omega]^{0,2}$$
is an $\Aut(\mn)$ equivariant map. 
Moreover, the orbit map $\tilde{p}:H^2(\mn)/\Aut(\mn)\lra E_\infty^{0,2}(\mn)/\Aut(\mn) $ is surjective. 
\end{Proposition}

\begin{Proof} The fact that $d(\mn^*)\subseteq \Lambda^2V_{k-1}$ implies that the map $p:H^2(\mn)\lra E_\infty^{0,2}$, $p\,([\omega])=[\omega]^{0,2}$ is well defined and surjective. Hence so is $\tilde{p}$.

 Notice that $p$ is injective if and only if $\dim H^2(\mn)=\dim E_\infty^{0,2}(\mn)$ and this situation occurs only when $\mn$ is a 2-step free nilpotent Lie algebra. 
\end{Proof}
\smallskip

In the next example we show that the quotient map $\tilde{p}$ is not always injective, even when $E_\infty^{0,2}\neq0$. 
\begin{Example} Let $\mn$ be the six dimensional nilpotent Lie algebra having non-zero Lie brackets
$$[e_1,e_2]=e_4,\quad [e_1,e_3]=e_5,\quad [e_1,e_4]=e_6. $$
Denote by $\{e^1,\ldots,e^6\}$ the dual basis of $\mn^*$. The followings are symplectic forms on $\mn$:
$$\omega_1=e^1\wedge e^6-e^2\wedge e^4+e^3\wedge e^5,\quad \omega_2=e^1\wedge e^6+e^2\wedge e^5+e^3\wedge e^4.$$
They verify $0\neq [\omega_1]^{0,2}=[\omega_2]^{0,2}=[e^1\wedge e^6]^{0,2}$. However through direct computations one can prove that the de Rham cohomology classes of $\omega_1$ and $\omega_2$ do not belong to the same $\Aut(\mn)$-orbit.
\end{Example}

\section{Classification of symplectic nilradicals.}

In this section we study the intermediate cohomology of the real nilpotent Lie algebras $\mn$ arising as nilradicals of minimal parabolic subalgebras of the real split forms of semisimple complex Lie algebras $\mgg$. In particular we show that $E_\infty^{0,2}(\mn)=0$ (in almost every case) when considering $\mgg$ to be a classical simple complex Lie algebra. According to the results in the previous section those nilpotent Lie algebras do not admit symplectic structures. Even more, we prove that if $E_\infty^{0,2}(\mn)\neq 0$ then $\mn$ does admit such structures.

To determine the intermediate cohomology group $E_\infty^{0,2}$ of those nilpotent Lie algebras, our main tool is the root decomposition of  semisimple Lie algebras. For this subject we give the book of Helgason \cite{HE} as a reference. The understanding of those systems allows the description of the filtration in Eq. (\ref{eq:eq53}) in terms of the root spaces. 
\bigskip

Let $\mg$ be a semisimple complex Lie algebra and let $\triangle$ be a root system of $\mg$. Then $\mg=\mh\oplus \bigoplus_{\alpha\in\triangle-\{0\}}\mg_\alpha$, where $\mh$ is a Cartan subalgebra of $\mg$. Denote as $\triangle^+$ the set of positive roots. The Lie algebra 
\begin{equation}\label{eq:complexnilp}
\mn=\bigoplus_{\alpha\in\triangle^+}\mg_\alpha\end{equation} is complex nilpotent. 

\begin{Remark} Given a complex nilpotent Lie algebra $\mn$, the filtration described in Eq. \prettyref{eq:eq53} and the induced spectral sequence are also canonically determined by $\mn$. Therefore each Lie algebra cohomology group with complex coefficients $H^i(\mn,\C)$ decomposes as in Eq. \prettyref{eq:eq8}.

In the particular case that $\mn$ is the nilradical in \prettyref{eq:complexnilp} of a Borel subalgebra of a complex semisimple Lie algebra, Kostant proved that $H^i(\mn,\C)$
is a direct sum of $T$-modules of dimension one (see \cite[Theorem 6.1]{KN2}) where $T$ is the diagonal subgroup of the semisimple Lie group. The action of $T$ on $\mn$ can be induced to $\Lambda^k\mn^*$ and commutes with the Lie algebra differential. 
As a consequence, the canonical filtration of $\mn^*$ is preserved by the $T$-action and the intermediate cohomology groups $E_\infty^{p,q}$ are $T$-modules. But $E_\infty^{p,q}$ is not irreducible in general. In particular, each (complex) intermediate cohomology group of $\mn$ is a sum of Kostant's one dimensional modules.
\end{Remark}

The Lie algebra $\mn$ in Eq. \prettyref{eq:complexnilp} admits a basis $\{X_\alpha\}_{\alpha\in\triangle^+}$ such that $\mg_\alpha=\C\,X_\alpha$ and the structure constants of $\mn$ in this basis are in $\R$. The object of study in this section is the real nilpotent Lie algebra having those real structure coefficients; we will also denote it as $\mn$. 

This real nilpotent Lie algebra $\mn$ is the nilradical of the minimal parabolic subalgebra of the split form corresponding to the semisimple Lie algebra $\mgg$.

We pursue the computation of the filtration of the Chevalley-Eilenberg complex of $\mn$  and the intermediate cohomology group $E_\infty^{0,2}(\mn)$.\vspace{0.2cm}

Denote as 
$\triangle_0=\{\alpha_1,\ldots,\alpha_r\}$ the subset of positive simple roots of $\mg$. Then for any positive root $\alpha$ there are non-negative integers $n_i$, $i=1,\ldots,r$ such that
$$\alpha=\sum_{i=1}^r n_i\alpha_i.$$
In this case we say that the {\em level} of the root is $\ell(\alpha)=\sum_{i=1}^rn_i$. Clearly the roots of level 1 are only the simple roots. There exists a unique positive root $\alpha_{\max}$ of maximal level, that is, such that $\ell(\alpha)\leq \ell(\alpha_{\max})$ for all $\alpha\in\triangle^+$. 

For each  $i\in\N$ define $L_i=\bigoplus _{\alpha:\ell(\alpha)=i} \R X_\alpha$ where $X_\alpha$ is as before. Then
$$\mn=\bigoplus_{j\geq 1} L_j\quad \mbox{ and }\quad[L_j,L_i]\subseteq L_{i+j}\,$$
and $\mn$ is an $\N$-graded Lie algebra.

Let $\alpha$ be a positive root of level $m$ and let $\gamma_\alpha\in\mn^*$ be  the dual element of $X_\alpha$. Since $[\mg_i,\mg_j]=\mg_m$ for all $i+j=m$, it holds 
\begin{equation}\label{eq:eq60}
d\gamma_\alpha\in \bigoplus_{i+j=m} L_i^*\wedge L_j^*.\end{equation} In particular $d\gamma_\alpha=0$ if and only if $\alpha$ is a simple root. This accounts into a description of the subspaces in (\ref{eq:eq53}) of $\mn^*$ as follows
$$V_0=0,\qquad V_j=\mbox{span}\{\gamma_\beta:\ell(\beta)\leq j\}=L_1^*\oplus\cdots\oplus L_j^*, \quad j=1,\ldots, k.$$
Notice that the nilpotency index of $\mn$ is $k=\ell(\alpha_{\max})$.

\smallskip

We proceed by making an insight into the space of closed 2-forms which we denote as $Z^2$. Such a form is an element of  $\Lambda^2\mn^*=\bigoplus_{1\leq i<j\leq k}L_i^*\wedge L_j^*$. In this context, the result of Benson and Gordon \cite[Lemma 2.8]{BE-GO} assures that $Z^2$ is contained in a strictly smaller subspace:
$$ Z^2\subseteq L_k^*\wedge L_1^* \oplus \bigoplus_{1\leq i<j\leq k-1}L_i^*\wedge L_j^*. $$ 

Therefore any $\omega\in Z^2$ can be written as a sum $\omega=\sigma+\tomega$ where 
\begin{equation}\label{eq:sigma}\sigma\in L_k^*\wedge L_1^*,\qquad \tomega\in \bigoplus_{1\leq i<j\leq k-1}L_i^*\wedge L_j^*\quad \mbox{and}\quad d\sigma=-d\tomega.\end{equation}

The vector space $L_k^*$ has dimension one and is spanned by $\gamma_{\alpha_{\max}}$. Moreover $L_1^*=V_1$ is spanned by the 1-forms $\gamma_{\alpha_i}$, $i=1,\ldots,r$ where $\alpha_1,\ldots,\alpha_r$  are the simple roots. Hence $\sigma= \gamma_{\alpha_{\max}}\wedge \eta$, with $\eta=\sum_{i=1}^r r_i \,\gamma_{\alpha_i}\,\in\,V_1$, $r_i\in\R$.

Since $d\eta=0$, $d\sigma=d \gamma_{\alpha_{\max}}\wedge \eta$ and by Eq. \prettyref{eq:eq60},  \begin{equation}\label{eq:eq63}
d\sigma \;\in\;L_{k-1}^*\wedge L_1^*\wedge L_1^* \;\oplus\; \bigoplus_{\stackrel{i+j=k}{1\leq i<j\leq k-2}}L_i^*\wedge L_j^*\;\wedge L_1^*.\end{equation}
In addition, \begin{equation}d\tomega\in\Lambda^3(L_{k-1}^*\oplus L_{k-2}^*\oplus\cdots\oplus L_1^*).\label{eq:dtomega} \end{equation}

The key here is to compare components of $d\sigma$ and $-d\omega$ in particular subspaces of those in Eqns. \prettyref{eq:eq63} and \prettyref{eq:dtomega}. For classical complex simple Lie algebras this idea allows us to prove that $Z^2\subset \Lambda^2V_{k-1}$ which by (\ref{eq6}) implies  $E_\infty^{0,2}(\mn)=0$. 

\begin{Remark} The structure of the semisimple Lie algebra $\mgg$ is independent of the root system. Therefore, the real nilpotent Lie algebra associated to a certain root system is isomorphic to the nilpotent Lie algebra arising from a different one. For this reason we can choose the root system of $\mgg$ that is more convenient for us to make calculations easier. \end{Remark}

The result we obtain is the following.



\begin{Lemma}\label{teo:teo32} Let $\mn$ be the nilradical of a minimal parabolic subalgebra of the real split Lie algebra corresponding to the complex classical simple Lie algebra $\mgg$.
\begin{enumerate}
\item If $\mgg=\mathfrak{sl}\,(n+1,\C)$ for some $n\geq 3$ then $E_\infty^{0,2}(\mn)=0$. 
\item If $\mgg=\mathfrak{so}\,(2n+1,\C)$ for some $n\geq 3$ then $E_\infty^{0,2}(\mn)=0$.
\item If $\mgg=\mathfrak{sp}\,(2n,\C)$ for some $n\geq 3$ then $E_\infty^{0,2}(\mn)=0$.
\item If $\mgg=\mathfrak{so}\,(2n,\C)$ for some $n\geq 4$ then $E_\infty^{0,2}(\mn)=0$.
\end{enumerate} 
\end{Lemma}

Using this Lemma, the full classification of the nilradicals admitting symplectic structures is obtained.

\begin{Theorem} \label{teo:teorema} Suppose $\mn$ is a nilradical of a minimal parabolic subalgebra of the real split Lie algebra corresponding to the complex classical simple Lie algebra $\mgg$. The Lie algebra of even dimension $\R^s\oplus \mn$, $s\geq 0$ admits symplectic structures if and only if $\mgg$ is one of the followings:
$$\mathfrak{sl}\,(2,\C), \quad\mathfrak{sl}\,(3,\C),\quad \mathfrak{so}\,(5,\C). $$ 
\end{Theorem}

\begin{Proof}
The nilradical $\mn$ of $\mathfrak{sl}\,(2,\C)$ is the abelian Lie algebra of dimension one and clearly $\R\oplus\mn$ admits symplectic structures. When $\mgg=\mathfrak{sl}\,(3,\C)$, then $\mn$ is the Heisenberg Lie algebra of dimension 3. It is well known that $\R\oplus \mn$ admits symplectic structures in this case also.

The nilradical corresponding to $\mathfrak{so}\,(5,\C)$ is the 4-dimensional 3-step nilpotent Lie algebra which can be endowed with a symplectic structure.

When $\mgg$ is not one of the Lie algebras above, Lemma \ref{teo:teo32} implies that its nilradical $\mn$ has $E_\infty^{0,2}(\mn)=0$. After Theorem \ref{teo:teo3}, those nilpotent Lie algebras are not symplectic.
\end{Proof}
\medskip

We proceed with the proof of Lemma \ref{teo:teo32}. This is  made using the canonical root systems known for classical simple Lie algebras. Moreover it is performed separately by cases because of the differences between those root systems. The order in which the cases are exposed is from the easiest to the most difficult one.

\smallskip
\begin{Proof} Part (1) of Theorem \ref{teo:teo32}. Let $\mg=\mathfrak{sl}(n+1,\C)$, $n\geq3$. If $n=3$, the Lie algebra $\mn$ is isomorphic to the Lie algebra of strictly upper triangular matrices $4\times 4$ for which can be easily shown that $E_\infty^{0,2}(\mn)=0$.

Suppose $n\geq 4$ and consider the Cartan subalgebra $\mh$ for which the positive roots are $e_i\pm e_j$, $1\leq i<j\leq n+1$. The set of simple roots is $\triangle_0=\{\alpha_i=e_i-e_{i+1}:i=1,\ldots,n\}$. Moreover the maximal root is $\alpha_{\max}=\sum_{i=1}^{n} \alpha_i$, hence the nilpotency index $k$ of $\mn$ is $k=n$. There are two different roots of level $n-1$, namely $\delta_1=\sum_{i=1}^{n-1}\alpha_i$ and $\delta_2=\sum_{i=2}^{n}\alpha_i$; in particular $\dim L_{k-1}=2$.

Let $\omega$ be a closed 2-form in $\mn^*$, then $\omega=\sigma+\tomega$ where $\sigma$ and $\tomega$ satisfy the conditions in \prettyref{eq:sigma}. The fact that $d\sigma=-d\tomega$ implies that the components of $d\sigma$ and $-d\tomega$ in the subspace $L_{k-1}^*\wedge L_1^*\wedge L_1^*$ are equal. So we compute both components.

As before, $\sigma=\gamma_{\alpha_{\max}}\wedge \eta$ where $\eta=\sum_{i=1}^{n} r_i\gamma_{\alpha_i}$, $r_i\in\R$ $i=1,\ldots,n$ and $d\sigma= d\gamma_{\alpha_{\max}}\wedge \eta$. 
Using Equation \prettyref{eq:eq60} and the fact that $\alpha_{\max}=\delta_1+\alpha_n=\delta_2+\alpha_1$ we obtain $d\gamma_{\alpha_{\max}}=a_1\, \gamma_{\delta_1}\wedge \gamma_{\alpha_{n}}+a_2\, \gamma_{\delta_2}\wedge \gamma_{\alpha_1} + \,\tau$ with $\tau\in\Lambda^2V_{k-2}=\Lambda^2(L_{k-2}^*\oplus\cdots\oplus L_1^*)$ and $a_1,a_2\in\R$ are both non-zero. This implies that the component of $d\sigma=d\gamma_{\alpha_{\max}}\wedge \eta$ in $ L_{k-1}^* \wedge L_1^*\wedge L_1^*$ is
\begin{equation}\label{eq:eq65}
  \sum_{i=1}^{n} a_1 r_i \, \gamma_{\delta_1}\wedge \gamma_{\alpha_{n}}\wedge \gamma_{\alpha_i}+\sum_{i=1}^n a_2 r_i \,\gamma_{\delta_2}\wedge \gamma_{\alpha_1} \wedge \gamma_{\alpha_i}
\end{equation}

To find the component of $d\tomega$ in the same subspace, write $\tomega=\omega_1+\tomega_1$ with
\begin{equation}\label{eq:eq62}
\omega_1\,\in\, L_{k-1}^*\wedge L_2^*
,\qquad \tomega_1\,\in\, L_{k-1}^*\wedge (\bigoplus_{\stackrel{j\leq k-1}{j\neq 2}} L_j^*\;) \oplus\, \Lambda^2(L_{k-2}^*\oplus\cdots\oplus L_1^*). 
\end{equation}
Hence by Eq. \prettyref{eq:eq60},
\begin{equation}\label{eq:eq64}
d\omega_1\,\in\, L_{k-1}^*\wedge L_1^*\wedge L_1^* \oplus \bigoplus_{i+j=k-1}L_i^*\wedge L_j^* \wedge L_2^*.
\end{equation}

By the same equation $d\tomega_1$ has no component in that subspace. Therefore the component of $d\sigma$ equals the component of $-d\omega_1$ in $L_{k-1}^*\wedge L_1^*\wedge L_1^*$.


Since $n\geq 4$ it is $L_{k-1}^*\neq L_2^*$ and the roots of level 2 are $\{\alpha_i+\alpha_{i+1}:i=1,\ldots,n-1\}$. Hence $L_{k-1}^*\wedge L_2^*$ has a basis of the form $\{\gamma_{\delta_1}\wedge \gamma_{\alpha_i+\alpha_{i+1}}, \gamma_{\delta_2}\wedge \gamma_{\alpha_i+\alpha_{i+1}}:i=1,\ldots,n-1\}$. Then $$\omega_1= \sum_{i=1}^{n-1} b_i\, \gamma_{\delta_1}\wedge \gamma_{\alpha_i+\alpha_{i+1}} + \sum_{i=1}^{n-1} c_i \,\gamma_{\delta_2}\wedge \gamma_{\alpha_i+\alpha_{i+1}},
\quad \mbox{ where } b_i\,c_i\,\in\R,\;\forall\,i.$$

Equation \prettyref{eq:eq60} implies that $d\gamma_{\alpha_i+\alpha_{i+1}}=\xi_i \gamma_{\alpha_i}\wedge \gamma_{\alpha_{i+1}}$ where $\xi_i\neq 0$, for all $i=1,\ldots,n-1$. 
Then
\begin{eqnarray}
-d\omega_1&=& \sum_{i=1}^{n-1} b_i \xi_i \gamma_{\delta_1}\wedge \gamma_{\alpha_i}\wedge\gamma_{\alpha_{i+1}} + \sum_{i=1}^{n-1} c_i\xi_i \gamma_{\delta_2}\wedge \gamma_{\alpha_i}\wedge\gamma_{\alpha_{i+1}}\nonumber\\
&&\quad
-\sum_{i=1}^n b_i d\gamma_{\delta_1}\wedge \gamma_{\alpha_i+\alpha_{i+1}} - \sum_{i=1}^n c_i d\gamma_{\delta_2}\wedge \gamma_{\alpha_i+\alpha_{i+1}}\nonumber.\end{eqnarray}

The component of $-d\omega_1$ in $L_{k-1}^*\wedge L_1^*\wedge L_1^*$ is
\begin{equation}\label{eq:eq66}
\sum_{i=1}^{n-1} b_i \xi_i \gamma_{\delta_1}\wedge \gamma_{\alpha_i}\wedge\gamma_{\alpha_{i+1}} + \sum_{i=1}^{n-1} c_i\xi_i \gamma_{\delta_2}\wedge \gamma_{\alpha_i}\wedge\gamma_{\alpha_{i+1}}.
\end{equation}
Indeed, the element $\;d\gamma_{\delta_i}$ belongs to $\Lambda^2V_{k-2}$ for $i=1,2$ because $\gamma_{\delta_i}\in L_{k-1}^*\subseteq V_{k-1}$. 

Formulas \prettyref{eq:eq65} and \prettyref{eq:eq66} give the components of $d\sigma$ and $-d\omega_1$ in $L_{k-1}^*\wedge L_1^*\wedge L_1^*$ which must be equal, that is,
$$
\sum_{i=1}^{n-1} b_i \xi_i \gamma_{\delta_1}\wedge \gamma_{\alpha_i}\wedge\gamma_{\alpha_{i+1}} + \sum_{i=1}^{n-1} c_i\xi_i \gamma_{\delta_2}\wedge \gamma_{\alpha_i}\wedge\gamma_{\alpha_{i+1}}
=\sum_{i=1}^{n-1} a_1 r_i\gamma_{\delta_1}\wedge \gamma_{\alpha_{n}}\wedge \gamma_{\alpha_i}+\sum_{i=2}^n a_2r_i \gamma_{\delta_2}\wedge \gamma_{\alpha_1} \wedge \gamma_{\alpha_i}.
$$

This expression implies that $r_{i}=0$ for all $i=1,\ldots,n$. Then $\eta=0$ implying $\sigma=0$. 
The conclusion here is that if $\omega=\sigma +\tomega$ is closed then $\sigma=0$, that is, $\omega\in \Lambda^2V_{k-1}$ and therefore $E_\infty^{0,2}(\mn)=0$ as we wanted to prove.\end{Proof}

\bigskip

\begin{Proof} Part (3) of Theorem \ref{teo:teo32}. Let $\mg=\mathfrak{sp}(2n,\C)$  for some $n\geq 3$. 
We prove the statement for $n\geq 4$.

Let $\mh$ be the Cartan subalgebra of $\mg$ associated to the set of positive roots $\triangle^+=\{e_i\pm e_j:1\leq i<j\leq n, \,2e_i,\,i=1,\ldots,n\}$; 
the subset of simple roots is $\triangle_0=\{\alpha_i:=e_i-e_{i+1}:i=1,\ldots,n-1,\alpha_n:=2e_n\}$. The maximal root is $\alpha_{\max}=\sum_{i=1}^{n-1} 2\alpha_i+\alpha_n$ from which we deduce that $\mn$ is $2n-1$-step nilpotent; set $k=2n-1$. Unlike the previous case, $\dim L_{k-1}=1$; the root of level $k-1$ is $\delta= \alpha_1+\sum_{i=2}^{n-1}2\alpha_i +\alpha_n$.
\medskip

Consider a closed 2-form $\omega\in\Lambda^2\mn^*$ with $\omega=\sigma +\tomega$ as in Equation (\ref{eq:sigma}). We are interested in the components of $d\sigma$ and $d\omega$ in the subspace $L_{k-1}^*\wedge L_1^*\wedge L_1^* \; \oplus\;L_{k-2}^*\wedge L_2^* \wedge L_1^*$ which must be opposite. Below we compute them both.

As before $\sigma=\gamma_{\alpha_{\max}}\wedge \eta$ where  $\eta$ is a linear combination of the 1-forms $\gamma_{\alpha_i}$, $i=1,\ldots,n$. Moreover $\alpha_{\max}=\delta +\alpha_1$. 
The roots $\rho=\alpha_1+\alpha_2+\sum_{i=3}^{n-1}2\alpha_i+\alpha_n$ and $\rho'=\sum_{i=2}^{n-1}2\alpha_i+\alpha_{n}$ are the only ones of level $k-2$ so
\begin{equation}\label{eq:eq67}d\gamma_{\alpha_{\max}}=a_1\gamma_{\alpha_1}\wedge\gamma_\delta+a_2\gamma_{\alpha_1+\alpha_2}\wedge\gamma_\rho +\tau,\qquad a_1\neq 0,a_2 \neq 0,\;  \mbox{ with } \tau\in\Lambda^2V_{k-3}.
\end{equation}
In fact, $\alpha_{\max}=\delta +\alpha_1$ and there do not exists any positive root $\beta$ such that $\beta+\rho'=\alpha_{\max}$. 


From Eq. \prettyref{eq:eq67}, 
$$d\sigma=d\gamma_{\alpha_{max}}\wedge \eta= a_1\gamma_{\alpha_1}\wedge\gamma_\delta\wedge \eta +a_2\gamma_{\alpha_1+\alpha_2}\wedge\gamma_\rho\wedge \eta +\tau\wedge \eta,\mbox{ being } \tau\wedge \eta\in \Lambda^3V_{k-3}$$
and its component in $L_{k-1}^* \wedge L_1^*\wedge L_1^*\oplus L_{k-2}^*\wedge L_2^*\wedge L_1^*$ is
\begin{equation}\label{eq:eq55}
a_1\gamma_{\alpha_1}\wedge\gamma_\delta\wedge \eta +a_2\gamma_{\alpha_1+\alpha_2}\wedge\gamma_\rho\wedge \eta. \end{equation}

To find the component of $d\tomega$ let $\omega_1$ and $\tomega_1$ be 2-forms such that $\tomega=\omega_1+\tomega_1$ where
$$\omega_1\in L_{k-1}^*\wedge L_2^*\, \oplus\, L_3^*\wedge L_{k-2}^*,\quad \mbox{and}\quad
\tomega_1\in L_{k-1}^*\wedge ( \bigoplus_{\stackrel{j\leq k-1}{j\neq 2}} L_j^*\;) \oplus L_{k-2}^*\wedge (\bigoplus_{\stackrel{j\leq k-1}{j\neq 3}} L_j^*\;) \,\oplus\, \Lambda^2(L_{k-3}^*\oplus\cdots\oplus L_1^*). 
$$
Using Eq. \prettyref{eq:eq60} for the differential of basic 1-forms one obtains
$$ d\omega_1\,\in\, L_{k-1}^*\wedge L_1^*\wedge L_1^* \; \oplus\;L_{k-2}^*\wedge L_2^* \wedge L_1^*.
$$
and $d\tomega_1$ has zero component in the same subspace. Hence the component of $d\tomega$ in $L_{k-1}^*\wedge L_1^*\wedge L_1^* \; \oplus\;L_{k-2}^*\wedge L_2^* \wedge L_1^*$ is the component of $d\omega_1$ in that subspace.


In this case, $\ell(\beta)=2$ if and only if $\beta=\alpha_i +\alpha_{i+1}$ for some $i=1,\ldots,n-1$. Thus $L_{k-1}^*\wedge L_2^*$ admits $\{\gamma_{\alpha_i+\alpha_{i+1}}\wedge \gamma_{\delta},\;i=1,\ldots, n-1\}$ as a basis. The roots of level three are $\alpha_i+\alpha_{i+1}+\alpha_{i+2}$, $i=1,\ldots,n-2$ and $2\alpha_{n-1}+\alpha_n$. 

Since $n\geq 4$, $3\neq k-2$ and $L_3^*\wedge L_{k-2}^*$ is spanned by $\{\gamma_{\alpha_i+\alpha_{i+1}+\alpha_{i+2}}\wedge \gamma_{\rho}, \gamma_{\alpha_i+\alpha_{i+1}+\alpha_{i+2}}\wedge \gamma_{\rho'}
:i=1,\ldots,n-2,\;
\gamma_{2\alpha_{n-1}+\alpha_n} \wedge\gamma_{\rho}
,\gamma_{2\alpha_{n-1}+\alpha_n} \wedge\gamma_{\rho'}
\}$.
The 2-form $\omega_1$ can be written as
\begin{eqnarray}
\omega_1&= &\sum_{i=1}^{n-1}c_i \gamma_{\alpha_i+\alpha_{i+1}}
\wedge \gamma_{\delta} + \sum_{i=1}^{n-2}d_i \gamma_{\alpha_i+\alpha_{i+1}+\alpha_{i+2}}\wedge \gamma_{\rho} + d_{n-1} \gamma_{2\alpha_{n-1}+\alpha_n} \wedge\gamma_{\rho}\nonumber\\
&& \qquad\qquad+ \sum_{i=1}^{n-2}e_i \gamma_{\alpha_i+\alpha_{i+1}+\alpha_{i+2}}\wedge \gamma_{\rho'} + e_{n-1} \gamma_{2\alpha_{n-1}+\alpha_n} \wedge\gamma_{\rho'},\quad c_i,d_i,e_i\in \R\nonumber.
\end{eqnarray}

From (\ref{eq:sigma}),
$$\hspace{-3cm}\begin{array}{lcl}
d\gamma_{\alpha_i+\alpha_{i+1}}&=&\xi_i \gamma_{\alpha_i}\wedge \gamma_{\alpha_{i+1}},\\
d\gamma_{\alpha_i+\alpha_{i+1}+\alpha_{i+2}}&=&s_i^1 \gamma_{\alpha_i+\alpha_{i+1}}\wedge\gamma_{\alpha_{i+2}}+ s_i^2\gamma_{\alpha_i}\wedge \gamma_{\alpha_{i+1}+\alpha_{i+2}}; \\
d\gamma_{2\alpha_{n-1}+\alpha_n}&=&s_{n-1}\gamma_{\alpha_n-1}\wedge\gamma_{\alpha_{n-1}+\alpha_n},
\end{array}$$
with $\xi_i,\,s_i^1,\,s_i^2,\,s_{n-1}$ all non-zero.

Putting this all together we reach:
\begin{eqnarray}
d\omega_1&=& \sum_{i=1}^{n-1}c_i\xi_i \gamma_{\alpha_i}\wedge \gamma_{\alpha_{i+1}}\wedge \gamma_{\delta} - \sum_{i=1}^{n-1}c_i \gamma_{\alpha_i+ \alpha_{i+1}}\wedge d\gamma_{\delta}\nonumber
\nonumber \\
&+&\sum_{i=1}^{n-2}d_i (s_i^1\gamma_{\alpha_i+\alpha_{i+1}}\wedge \gamma_{ \alpha_{i+2}} +s_i^2 \gamma_{\alpha_i}\wedge\gamma_{\alpha_{i+1}+\alpha_{i+2}})\wedge \gamma_{\rho} \nonumber\\ 
&-&\sum_{i=1}^{n-2}d_i\gamma_{\alpha_i+\alpha_{i+1}+\alpha_{i+2}}\wedge d\gamma_\rho + d_{n-1} s_{n-1}\gamma_{\alpha_{n-1}}\wedge \gamma_{\alpha_{n-1}+\alpha_n}  \wedge\gamma_{\rho} -d_{n-1}\gamma_{2\alpha_{n-1}+\alpha_n}\wedge d\gamma_{\rho}+
 \nonumber\\
 &-&\sum_{i=1}^{n-2}e_i \gamma_{\alpha_i+\alpha_{i+1}+\alpha_{i+2}}\wedge d\gamma_{\rho'} 
+ \sum_{i=1}^{n-2}e_i (s_i^1\gamma_{\alpha_i+\alpha_{i+1}}\wedge \gamma_{\alpha_{i+2}} +s_i^2 \gamma_{\alpha_i}\wedge\gamma_{\alpha_{i+1}+\alpha_{i+2}})\wedge \gamma_{\rho'} \nonumber\\
&-& e_{n-1} s_{n-1}\gamma_{\alpha_{n-1}}\wedge\gamma_{\alpha_{n-1}+\alpha_n}\wedge\gamma_{\rho'} -e_{n-1} \gamma_{2\alpha_{n-1}+\alpha_n}\wedge d\gamma_{\rho'}.\nonumber
\end{eqnarray}

Since $\ell(\rho)=\ell(\rho')=k-2$, $d\gamma_\rho,\,d\gamma_{\rho'}\in\Lambda^2V_{k-3}$. In addition,  $\delta=\alpha_2+\rho=\alpha_1+\rho'$ implies
\begin{eqnarray}
d\gamma_{\delta}&=&b_1\gamma_{\alpha_2} \wedge\gamma_\rho +b_2\gamma_{\alpha_1}\wedge\gamma_{\rho'} +\tau',\qquad b_1\neq 0,b_2 \neq 0,\; \mbox{ with }\tau'\in\Lambda^2 V_{k-3}. \nonumber
\end{eqnarray}

Therefore, the component of $-d\omega_1$ in $ L_{k-1}^* \wedge L_1^*\wedge L_1^*\oplus L_{k-2}^*\wedge L_2^*\wedge L_1^*$ is
\begin{eqnarray}
&&\hspace{-1cm}-\sum_{i=1}^{n-1}c_i\xi_i \gamma_{\alpha_i}\wedge \gamma_{\alpha_{i+1}}\wedge \gamma_{\delta} + \sum_{i=1}^{n-1}c_i \gamma_{\alpha_i+ \alpha_{i+1}}\wedge (b_1\gamma_{\alpha_2} \wedge\gamma_\rho +b_2\gamma_{\alpha_1}\wedge\gamma_{\rho'})-
 \label{eq:eq56} \\
&-& \sum_{i=1}^{n-2}d_i (s_i^1\gamma_{\alpha_i+\alpha_{i+1}}\wedge \gamma_{ \alpha_{i+2}} +s_i^2 \gamma_{\alpha_i}\wedge\gamma_{\alpha_{i+1}+\alpha_{i+2}})\wedge \gamma_{\rho} - d_{n-1} s_{n-1}\gamma_{\alpha_{n-1}}\wedge \gamma_{\alpha_{n-1}+\alpha_n}  \wedge\gamma_{\rho}-\nonumber\\
&-&\sum_{i=1}^{n-2}e_i (s_i^1\gamma_{\alpha_i+\alpha_{i+1}}\wedge \gamma_{\alpha_{i+2}} +s_i^2 \gamma_{\alpha_i}\wedge\gamma_{\alpha_{i+1}+\alpha_{i+2}})\wedge \gamma_{\rho'} + e_{n-1} s_{n-1}\gamma_{\alpha_{n-1}}\wedge\gamma_{\alpha_{n-1}+\alpha_n}\wedge\gamma_{\rho'} .\nonumber
\end{eqnarray}

The components of $d\sigma$ and $-d\tomega_1$ expressed in Eqns. \prettyref{eq:eq55} and \prettyref{eq:eq56} respectively, coincide. Following usual computations we get that
$c_i$, $e_j$ are zero for all $i=1,\ldots,n-1$ and $j=1,\ldots,n-2$ which simplify Eq. \prettyref{eq:eq56}. After some other simplifications we obtain once again that $\sigma=0$ and then 
any closed 2-form in $\mn$ belongs to $\Lambda^2 V_{k-1}$ which is equivalent to $E_\infty^{0,2}(\mn)=0$.

When $n=3$ the Lie algebra $\mn$ has dimension 9 and the proof is made by direct computations using a mathematical software.
\end{Proof}
\bigskip

\begin{Proof} Part (2) of Theorem \ref{teo:teo32}. Suppose $\mg=\mathfrak{so}(2n+1,\C)$ for some $n\geq 3$.

We perform the proof for $n\geq 5$. Consider the Cartan Lie subalgebra such that the set of positive roots is $\triangle^+=\{e_i\pm e_j:\,1\leq i<j\leq n, \,e_i,\,i=1,\ldots,n\}$. Then $\dim \mn= n^2$. The simple roots are $ \{\alpha_i:=e_i-e_{i+1}:i=1,\ldots,n-1,\alpha_n:=e_n\}$ and the maximal root $\alpha_{\max}$ is $e_1+e_2=\alpha_1+2\sum_{i=2}^{n} \alpha_i$. As a consequence $\mn$ is $2n-1$-step nilpotent; set $k=2n-1$. Here we also have $\dim L_{k-1}=1$; the root of level $k-1$ is now $\delta= \alpha_1+\alpha_2+\sum_{i=2}^{n-1}2\alpha_i $.

 Let $\omega\in\Lambda^2\mn^*$ be a closed form and let $\sigma$ and $\tomega$ be as in \prettyref{eq:sigma}. In this case we study the components of $d\sigma$ and $d\tomega$ in $L_{k-1}^*\wedge L_1^*\wedge L_1^* \; \oplus\;L_{k-2}^*\wedge L_2^* \wedge L_1^*\oplus L_{k-3}^*\wedge L_3^*\wedge L_1^*$ which must be opposite since $0=d\omega=d\sigma+d\tomega$.
 \smallskip
 
 Recall that $\sigma=\gamma_{\alpha_{\max}}\wedge \eta$ with $\eta$ a linear combination of the 1-forms $\gamma_{\alpha_i}$, $i=1,\ldots,n$. There are two roots of level $k-2$, namely $\rho=\alpha_1+\alpha_2+\alpha_3+\sum_{i=4}^{n}2\alpha_i$ and $\rho'=\alpha_2+\sum_{i=3}^{n}2\alpha_i$. The roots of level $k-3$ are $\theta=\alpha_1+\alpha_2+\alpha_3+\alpha_4+2\sum_{i=5}^n\alpha_i$ and $\theta'=\alpha_2+\alpha_3+2\sum_{i=4}^n\alpha_i$.

Notice that 
$$\alpha_{\max}=\delta +\alpha_2,\quad \alpha_{\max}=\rho+(\alpha_2+\alpha_3)=\rho'+(\alpha_1+\alpha_2)\quad  \mbox{ and }$$
$$\alpha_{\max}=\theta +(\alpha_2+\alpha_3+\alpha_4)=\theta'+(\alpha_1+\alpha_2+\alpha_3).$$ This implies
\begin{eqnarray}
d\gamma_{\alpha_{\max}}&=&a_1\,\gamma_{\alpha_2}\wedge\gamma_\delta+a_2 \,\gamma_{\alpha_2+\alpha_3}\wedge \gamma_\rho+a_3\,\gamma_{\alpha_1+\alpha_2}\wedge\gamma_{\rho'}\label{eq:eq85}\\
&&\quad+a_4\,\gamma_{\alpha_2+\alpha_3+\alpha_4}\wedge \gamma_\theta+a_5\,\gamma_{\alpha_1+\alpha_2+\alpha_3}\wedge \gamma_{\theta'} +\tau,\qquad a_i\neq 0,\,\forall i,\;\mbox{ and } \tau\in\Lambda^2V_{k-4}.\nonumber
\end{eqnarray}

Then $d\sigma=d\gamma_{\alpha_{\max}}\wedge \eta$ is
\begin{eqnarray}
d\sigma&=& a_1\,\gamma_{\alpha_2}\wedge\gamma_\delta\wedge\eta+a_2 \,\gamma_{\alpha_2+\alpha_3}\wedge \gamma_\rho\wedge\eta+a_3\,\gamma_{\alpha_1+\alpha_2}\wedge\gamma_{\rho'}\wedge\eta\nonumber\\
&&\quad\qquad+a_4\,\gamma_{\alpha_2+\alpha_3+\alpha_4}\wedge \gamma_\theta\wedge\eta+a_5\,\gamma_{\alpha_1+\alpha_2+\alpha_3}\wedge \gamma_{\theta'} \wedge\eta +\tilde{\tau},\qquad \;\mbox{ with }  \tilde{\tau}\in\Lambda^2V_{k-4}\nonumber.\end{eqnarray}

The component of $d\sigma$ in $L_{k-1}^*\wedge L_1^*\wedge L_1^* \; \oplus\;L_{k-2}^*\wedge L_2^* \wedge L_1^*\oplus L_{k-3}^*\wedge L_3^*\wedge L_1^*$ is
\begin{eqnarray}\hspace{-3cm}&&a_1\,\gamma_{\alpha_2}\wedge\gamma_\delta\wedge\eta+a_2 \,\gamma_{\alpha_2+\alpha_3}\wedge \gamma_\rho\wedge\eta+a_3\,\gamma_{\alpha_1+\alpha_2}\wedge\gamma_{\rho'}\wedge\eta\nonumber\\ &&\hspace{2cm}\qquad+a_4\,\gamma_{\alpha_2+\alpha_3+\alpha_4}\wedge \gamma_\theta\wedge\eta+a_5\,\gamma_{\alpha_1+\alpha_2+\alpha_3}\wedge \gamma_{\theta'} \wedge\eta .\label{eq:eq86}
\end{eqnarray}

To compute the component of $d\omega$ in the same subspace, take $\omega_1$ and $\tomega_1$ such that $\tomega=\omega_1+\tomega_1$ where
$$\omega_1\in  L_2^*\wedge L_{k-1}^*\, \oplus\, L_3^*\wedge L_{k-2}^*\,\oplus \,L_4^*\wedge L_{k-3}\qquad \mbox{and}$$
$$
\tomega_1\in L_{k-1}^*\wedge ( \bigoplus_{\stackrel{j\leq k-1}{j\neq 2}} L_j^*\;) \oplus L_{k-2}^*\wedge (\bigoplus_{\stackrel{j\leq k-1}{j\neq 3}} L_j^*\;)\oplus L_{k-3}^*\wedge ( \bigoplus_{\stackrel{j\leq k-3}{j\neq 4}} L_k^*\;) \,\oplus\, \Lambda^2(L_{k-4}^*\oplus\cdots\oplus L_1^*). 
$$

The 3-form $d\tomega_1$ has zero component in $L_{k-1}^*\wedge L_1^*\wedge L_1^* \; \oplus\;L_{k-2}^*\wedge L_2^* \wedge L_1^*\oplus L_{k-3}^*\wedge L_3^*\wedge L_1^*$. Therefore the component of $d\tomega$ in that subspace is that one of $d\omega_1$ and to compute it we will make use of the following formulas obtained from Eq. \prettyref{eq:eq60}:
\hspace{-4cm}\begin{eqnarray}
d\gamma_\delta &=& b_1\,\gamma_{\alpha_3}\wedge \gamma_\rho+b_2\,\gamma_{\alpha_1}\wedge \gamma_{\rho'}+ b_3\,\gamma_\theta\wedge \gamma_{\alpha_3+\alpha_4}+\tau',\nonumber\\
d\gamma_\rho&=&\nu_1\, \gamma_{\alpha_1}\wedge \gamma_{\theta'
}+\nu_2\,\gamma_{\alpha_4}\wedge \gamma_\theta +\tau'',\label{eq:eq87}\\
d\gamma_{\rho'}&=&\mu_1\, \gamma_{\alpha_3}\wedge \gamma_{\theta'
}+\tau''',\qquad\qquad \mbox{ con } \tau',\tau'',\tau''\in\Lambda^2V_{k-4}\nonumber 
\end{eqnarray}
The difference between $d\gamma_\rho$ and $d\gamma_{\rho'}$ is due to the lack of simple roots $\beta$ verifying $\theta+\beta=\rho'$; in opposite to $\rho$ which verifies $\rho=\theta+\alpha_4$. 

Notice that $\ell(\theta)=\ell(\theta')=k-3$, hence $\gamma_\theta$ and $\gamma_{\theta'}$ are in $V_{k-3}$;
$d\gamma_\theta$ and $d\gamma_{\theta'}$ are elements in $\Lambda^2V_{k-4}=\Lambda^2(L_{k-4}^*\oplus\cdots\oplus L_1^*)$.
\smallskip

The roots of level two are $\alpha_i +\alpha_{i+1}$, $i=1,\ldots,n-1$, which gives the following basis of  $L_{k-1}^*\wedge L_2^*$:  $\{\gamma_{\alpha_i+\alpha_{i+1}}\wedge \gamma_{\delta},\;i=1,\ldots, n-1\}$. In addition, the roots of level three are $\alpha_i+\alpha_{i+1}+\alpha_{i+2}$, $i=1,\ldots,n-2$ and $\alpha_{n-1}+2\alpha_n$. Therefore $L_3^*\wedge L_{k-2}^*$ is spanned by $\{\gamma_{\alpha_i+\alpha_{i+1}+\alpha_{i+2}}\wedge \gamma_{\rho}, \gamma_{\alpha_i+\alpha_{i+1}+\alpha_{i+2}}\wedge \gamma_{\rho'}
:i=1,\ldots,n-2,\;
\gamma_{\alpha_{n-1}+2\alpha_n} \wedge\gamma_{\rho}
,\gamma_{\alpha_{n-1}+2\alpha_n} \wedge\gamma_{\rho'}
\}$.
Finally, the roots of level four are $\alpha_i+\alpha_{i+1}+\alpha_{i+2}+\alpha_{i+3}$, $i=1,\ldots,n-3$, and $\alpha_{n-2}+\alpha_{n-1}+2\alpha_n$. Since $n\geq 5$, $k-3\neq 4$ and   $\{\gamma_{\alpha_i+\alpha_{i+1}+\alpha_{i+2}+\alpha_{i+3}}\wedge\gamma_\theta,\,\gamma_{\alpha_i+\alpha_{i+1}+\alpha_{i+2}+\alpha_{i+3}}\wedge\gamma_{\theta'},\,i=1,\ldots,n-3,\,\gamma_{\alpha_{n-2}+\alpha_{n-1}+2\alpha_n}\wedge \gamma_\theta,\,\gamma_{\alpha_{n-2}+\alpha_{n-1}+2\alpha_n}\wedge \gamma_{\theta'}\}$ is a basis of $L_{k-3}^*\wedge L_4^*$. 

To make computations easier, write $\omega_1\in L_2^*\wedge L_{k-1}^*\, \oplus\, L_3^*\wedge L_{k-2}^*\,\oplus \,L_4^*\wedge L_{k-3}^*$ as $\omega_1=\omega_1^a+\omega_1^b+\omega_1^c$ where for some coefficients $c_i,d_i,e_if_i,g_i\in \R$, it holds
\begin{eqnarray}
\omega_1^a&=& \sum_{i=1}^{n-1} c_i\,\gamma_{\alpha_i+\alpha_{i+1}}\wedge \gamma_{\delta}\; \in \,L_{k-1}^*\wedge L_2^*,\nonumber\\
\omega_1^b&=&\sum_{i=1}^{n-2}d_i\,\gamma_{\alpha_i+\alpha_{i+1}+\alpha_{i+2}}\wedge \gamma_{\rho}+d_{n-1}\,\gamma_{\alpha_{n-1}+2\alpha_n} \wedge\gamma_{\rho}+\nonumber\\
&&\quad+\sum_{i=1}^{n-2}e_i\,\gamma_{\alpha_i+\alpha_{i+1}+\alpha_{i+2}}\wedge \gamma_{\rho'}+ e_{n-1}\,\gamma_{\alpha_{n-1}+2\alpha_n} \wedge\gamma_{\rho'} \in\, L_3^*\wedge L_{k-2}^*,\nonumber\\
\omega_1^c&=&\sum_{i=1}^{n-3} f_i\,\gamma_{\alpha_i+\alpha_{i+1}+\alpha_{i+2}+\alpha_{i+3}}\wedge\gamma_\theta +f_{n-2}\,\gamma_{\alpha_{n-2}+\alpha_{n-1}+2\alpha_n}\wedge \gamma_\theta+\nonumber \\
&&\quad +\sum_{i=1}^{n-3} g_i\,\gamma_{\alpha_i+\alpha_{i+1}+\alpha_{i+2}+\alpha_{i+3}}\wedge\gamma_{\theta'} +g_{n-2}\,\gamma_{\alpha_{n-2}+\alpha_{n-1}+2\alpha_n}\wedge \gamma_{\theta'}\;\in\,L_4^*\wedge L_{k-3}^*.\nonumber
\end{eqnarray}

The differential of the basic elements of $L_2^*$, $L_3^*$ y $L_4^*$ are
$$\begin{array}{lcl}
d\gamma_{\alpha_i+\alpha_{i+1}}&=&\xi_i\, \gamma_{\alpha_i}\wedge \gamma_{\alpha_{i+1}},\\
d\gamma_{\alpha_i+\alpha_{i+1}+\alpha_{i+2}}&=&s_i^1 \gamma_{\alpha_i+\alpha_{i+1}}\wedge\gamma_{\alpha_{i+2}}+ s_i^2\gamma_{\alpha_i}\wedge \gamma_{\alpha_{i+1}+\alpha_{i+2}}, \\
d\gamma_{\alpha_{n-1}+2\alpha_n}&=&s_{n-1}\gamma_{\alpha_{n-1}+\alpha_n}\wedge\gamma_{\alpha_n},\\
d\gamma_{\alpha_i+\alpha_{i+1}+\alpha_{i+2}+\alpha_{i+3}}&=& t_i^1\; \gamma_{\alpha_i}\wedge \gamma_{\alpha_{i+1}+\alpha_{i+2}+\alpha_{i+3}}+t_i^2 \;\gamma_{\alpha_{i}+\alpha_{i+1}}\wedge \gamma_{\alpha_{i+2}+\alpha_{i+3}}+t_i^3\; \gamma_{\alpha_i+\alpha_{i+1}+\alpha_{i+2}}\wedge \gamma_{\alpha_{i+3}},\\
d\gamma_{\alpha_{n-2}+\alpha_{n-1}+2\alpha_n}&=& t_{n-2}^1 \,\gamma_{\alpha_{n-2}}\wedge \gamma_{\alpha_{n-1}+2\alpha_n}+t_{n-2}^2 \,\gamma_{ \alpha_{n-2}+\alpha_{n-1}+\alpha_{n}}\wedge \gamma_{\alpha_n}, 
\end{array}$$
where $\xi_i,\,s_i^j,s_{n-1},t_i^j,$ are all non-zero.

The differential of  $\omega_1^a$ is by \prettyref{eq:eq87}
$$d\omega_1^a= \sum_{i=1}^{n-1} c_i\,d\gamma_{\alpha_i+\alpha_{i+1}}\wedge \gamma_{\delta}-\sum_{i=1}^{n-1} c_i\,\gamma_{\alpha_i+\alpha_{i+1}}\wedge d\gamma_{\delta}= \sum_{i=1}^{n-1} c_i\,\xi_i\, \gamma_{\alpha_i}\wedge \gamma_{\alpha_{i+1}}\wedge \gamma_{\delta}-$$$$\sum_{i=1}^{n-1} c_i\,\gamma_{\alpha_i+\alpha_{i+1}}\wedge \left( b_1\,\gamma_{\alpha_3}\wedge \gamma_\rho+b_2\,\gamma_{\alpha_1}\wedge \gamma_{\rho'}+ b_3\,\gamma_\theta\wedge \gamma_{\alpha_3+\alpha_4}+\tau'\right).$$

Then the component of $d\omega_1^a$ in $L_{k-1}^*\wedge L_1^*\wedge L_1^* \; \oplus\;L_{k-2}^*\wedge L_2^* \wedge L_1^*\oplus L_{k-3}^*\wedge L_3^*\wedge L_1^*$ is
\begin{equation}\label{eq:eq88}\sum_{i=1}^{n-1} c_i\,\xi_i\, \gamma_{\alpha_i}\wedge \gamma_{\alpha_{i+1}}\wedge \gamma_{\delta}-
\sum_{i=1}^{n-1}  b_1\,c_i\,\gamma_{\alpha_i+\alpha_{i+1}}\wedge\gamma_{\alpha_3}\wedge \gamma_\rho+b_2\,c_i\,\gamma_{\alpha_i+\alpha_{i+1}}\wedge\gamma_{\alpha_1}\wedge \gamma_{\rho'},
\end{equation}
since $\tau'\in\Lambda^2V_{k-4}=\Lambda^2(L_{k-4}^*\oplus \cdots\oplus L_1^*)$ and $\gamma_{\alpha_i+\alpha_{i+1}}\wedge \gamma_\theta\wedge \gamma_{\alpha_3+\alpha_4}\in L_2^*\wedge L_2^*\wedge L_{k-3}^*$.

In a similar way we compute the components of $d\omega_1^b$ and $d\omega_1^c$ in $L_{k-1}^*\wedge L_1^*\wedge L_1^* \; \oplus\;L_{k-2}^*\wedge L_2^* \wedge L_1^*\oplus L_{k-3}^*\wedge L_3^*\wedge L_1^*$ which are:

\noindent - component of $d\omega_1^b$:
\begin{eqnarray}
&&\sum_{i=1}^{n-2}\left( d_is_i^1\,\gamma_{\alpha_i+\alpha_{i+1}}\wedge \gamma_{\alpha_{i+2}}\wedge \gamma_\rho+
 s_i^2d_i\,\gamma_{\alpha_i}\wedge \gamma_{\alpha_{i+1}+\alpha_{i+2}}\wedge\gamma_\rho\right)+
d_{n-1}s_{n-1}\,\gamma_{\alpha_{n-1}+\alpha_n}\wedge\gamma_{\alpha_n}\wedge \gamma_\rho+\nonumber \\
&&-\sum_{i=1}^{n-2}\left( d_i\nu_1\, \gamma_{\alpha_i+\alpha_{i+1}+\alpha_{i+2}}\wedge\gamma_{\alpha_1}\wedge\gamma_{\theta'}+ d_i\nu_2\, \gamma_{\alpha_i+\alpha_{i+1}+\alpha_{i+2}}\wedge\gamma_{\alpha_4}\wedge\gamma_{\theta}\right)-d_{n-1}\nu_1\, \gamma_{\alpha_{n-1}+2\alpha_n}\wedge \gamma_{\alpha_1}\wedge\gamma_{\theta'}+\nonumber\\ 
&&-\nu_2d_{n-1}\gamma_{\alpha_{n-1}+2\alpha_n}\wedge \gamma_{\alpha_4}\wedge \gamma_\theta+\sum_{i=1}^{n-2} d_is_i^1\,\gamma_{\alpha_i+\alpha_{i+1}}\wedge \gamma_{\alpha_{i+2}}\wedge\gamma_{\rho'}+
\sum_{i=1}^{n-2} d_is_i^2\,\gamma_{\alpha_i}\wedge \gamma_{\alpha_{i+1}+\alpha_{i+2}}\wedge\gamma_{\rho'}\label{eq:eq89}\\
&& d_{n-1}s_{n-1}\, \gamma_{\alpha_{n-1}+\alpha_n}\wedge \gamma_{\alpha_n}\wedge \gamma_{\rho'}-\sum_{i=1}^{n-2}d_i\mu_1\,\gamma_{\alpha_i+\alpha_{i+1}+\alpha_{i+2}}\wedge \gamma_{\alpha_3}\wedge\gamma_{\theta'}-d_{n-1}\mu_1\,\gamma_{\alpha_{n-1}+2\alpha_n}\wedge \gamma_{\alpha_3}\wedge \gamma_{\theta'};\nonumber
\end{eqnarray}
notice that, in fact, the elements of \prettyref{eq:eq89} belong to $L_{k-2}^*\wedge L_2^* \wedge L_1^*\oplus L_{k-3}^*\wedge L_3^*\wedge L_1^*$.

\noindent - component of $d\omega_1^c$ :
\begin{eqnarray}
&&\sum_{i=1}^{n-3}f_i\,\left(t_i^1\,\gamma_{\alpha_i}\wedge \gamma_{\alpha_{i+1}+\alpha_{i+2}+\alpha_{i+3}}+t_i^3\gamma_{\alpha_i+\alpha_{i+1}+\alpha_{i+2}}\wedge\gamma_{\alpha_{i+3}} \right)\wedge\gamma_\theta  +
f_{n-2}( t_{n-2}^1\gamma_{\alpha_{n-2}}\wedge\gamma_{\alpha_{n-1}+2\alpha_n}+\nonumber\\
&&+t_{n-2}^2\gamma_{\alpha_{n-2}+\alpha_{n-1}+\alpha_n}\wedge\gamma_{\alpha_n})\wedge \gamma_\theta +\sum_{i=1}^{n-3}g_i\,\left(t_i^1\,\gamma_{\alpha_i}\wedge \gamma_{\alpha_{i+1}+\alpha_{i+2}+\alpha_{i+3}}+t_i^3\gamma_{\alpha_i+\alpha_{i+1}+\alpha_{i+2}}\wedge\gamma_{\alpha_{i+3}} \right)\wedge\gamma_{\theta'}  \nonumber\\
&&+
g_{n-2}( t_{n-2}^1\gamma_{\alpha_{n-2}}\wedge\gamma_{\alpha_{n-1}+2\alpha_n} +t_{n-2}^2\gamma_{\alpha_{n-2}+\alpha_{n-1}+\alpha_n}\wedge \gamma_{\alpha_n})\wedge \gamma_{\theta'} .
\label{eq:eq90}
\end{eqnarray}
In this case, the elements of  \prettyref{eq:eq90} are in $L_{k-3}^*\wedge L_3^*\wedge L_1^*$.
 
The component of  $-d\omega_1$ in $L_{k-1}^*\wedge L_1^*\wedge L_1^* \; \oplus\;L_{k-2}^*\wedge L_2^* \wedge L_1^*\oplus L_{k-3}^*\wedge L_3^*\wedge L_1^*$ is obtained from Eqns. \prettyref{eq:eq88}, \prettyref{eq:eq89} and \prettyref{eq:eq90} and, at the same time, it coincides with 3-form in \prettyref{eq:eq86}.

The part in $L_{k-1}^*\wedge L_1^*\wedge L_1^*$ of \prettyref{eq:eq86} and that of $-d\omega_1$ are equal, that is,
$$\sum_{i=1}^n a_1\,r_i\,\gamma_{\alpha_2}\wedge\gamma_\delta\wedge\gamma_{\alpha_i} =\sum_{i=1}^{n-1} c_i\,\xi_i\, \gamma_{\alpha_i}\wedge \gamma_{\alpha_{i+1}}\wedge \gamma_{\delta}. $$
This imply $r_i=0$ for all $4\leq i\leq n$. Putting this in \prettyref{eq:eq86} and looking at the part in $ L_{k-3}^*\wedge L_3^*\wedge L_1^*$  of $-d\omega_1$ we obtain
\begin{eqnarray}
&&\sum_{i=1}^3\eta_i\,(a_4\,\gamma_{\alpha_2+\alpha_3+\alpha_4}\wedge \gamma_\theta+a_5\,\gamma_{\alpha_1+\alpha_2+\alpha_3}\wedge \gamma_{\theta'})\wedge\gamma_{\alpha_i}=\nonumber\\
&&\sum_{i=1}^{n-3}f_i\,\left(t_i^1\,\gamma_{\alpha_i}\wedge \gamma_{\alpha_{i+1}+\alpha_{i+2}+\alpha_{i+3}}+t_i^3\gamma_{\alpha_i+\alpha_{i+1}+\alpha_{i+2}}\wedge\gamma_{\alpha_{i+3}} \right)\wedge\gamma_\theta  +
f_{n-2}( t_{n-2}^1\gamma_{\alpha_{n-2}}\wedge\gamma_{\alpha_{n-1}+2\alpha_n}+\nonumber\\
&&+t_{n-2}^2\gamma_{\alpha_{n-2}+\alpha_{n-1}+\alpha_n}\wedge\gamma_{\alpha_n})\wedge \gamma_\theta +\sum_{i=1}^{n-3}g_i\,\left(t_i^1\,\gamma_{\alpha_i}\wedge \gamma_{\alpha_{i+1}+\alpha_{i+2}+\alpha_{i+3}}+t_i^3\gamma_{\alpha_i+\alpha_{i+1}+\alpha_{i+2}}\wedge\gamma_{\alpha_{i+3}} \right)\wedge\gamma_{\theta'}  \nonumber\\
&&+
g_{n-2}( t_{n-2}^1\gamma_{\alpha_{n-2}}\wedge\gamma_{\alpha_{n-1}+2\alpha_n} +t_{n-2}^2\gamma_{\alpha_{n-2}+\alpha_{n-1}+\alpha_n}\wedge \gamma_{\alpha_n})\wedge \gamma_{\theta'}\\
&&-\sum_{i=1}^{n-2}\left( d_i\nu_1\, \gamma_{\alpha_i+\alpha_{i+1}+\alpha_{i+2}}\wedge\gamma_{\alpha_1}\wedge\gamma_{\theta'}+ d_i\nu_2\, \gamma_{\alpha_i+\alpha_{i+1}+\alpha_{i+2}}\wedge\gamma_{\alpha_4}\wedge\gamma_{\theta}\right)-d_{n-1}\nu_1\, \gamma_{\alpha_{n-1}+2\alpha_n}\wedge \gamma_{\alpha_1}\wedge\gamma_{\theta'}+\nonumber\\ 
&&-\nu_2d_{n-1}\gamma_{\alpha_{n-1}+2\alpha_n}\wedge \gamma_{\alpha_4}\wedge \gamma_\theta +\sum_{i=1}^{n-2}d_i\mu_1\,\gamma_{\alpha_i+\alpha_{i+1}+\alpha_{i+2}}\wedge \gamma_{\alpha_3}\wedge\gamma_{\theta'}-d_{n-1}\mu_1\,\gamma_{\alpha_{n-1}+2\alpha_n}\wedge \gamma_{\alpha_3}\wedge \gamma_{\theta'}\nonumber\end{eqnarray}

Being careful and comparing term by term we deduce that $r_2=r_3=0$ and $f_i=d_i=0$ for all $i\geq 2$. Comparing one more time we reach $r_1=0$ and therefore $\sigma=0$. As before, this implies $E_\infty^{0,2}(\mn)=0$. 

If $n=3$ or $n=4$ then the intermediate cohomology groups $E_{\infty}^{0,2}(\mn)$ were proved to be zero throughout a computational program. 
\end{Proof}

\bigskip

\begin{Proof} Part (4) of Theorem \ref{teo:teo32}. The family of simple Lie algebras $\mathfrak{so}(2n,\C)$ is defined for $n\geq4$. 

Let $\mh$ be the Cartan subalgebra for which corresponds the following set of positive roots $\triangle^+=\{e_i\pm e_j:\,1\leq i<j\leq n\}$. Hence $\dim \mn= n(n-1)$.
The set of simple roots is $\triangle_0=\{\alpha_i:=e_i-e_{i+1}:i=1,\ldots,n-1,\alpha_n:=e_{n-1}+e_n\}$ and the maximal root $\alpha_{\max}$ is $e_1+e_2$ and can be obtained as $\alpha_{\max}=\alpha_1+2\sum_{i=2}^{n-2} \alpha_i+\alpha_{n-1}+\alpha_n$. Then $\mn$ is $2n-3$-step nilpotent. Define $k=2n-3$. As in the previous case $\dim L_{k-1}=1$ and the root of level $k-1$ is $\delta= \alpha_1+\alpha_2+\sum_{i=2}^{n-2}2\alpha_i+\alpha_{n-1}+\alpha_n$.
\smallskip

If $n\geq 6$, then the proof of the previous case applies. Actually, the root system corresponding to $\mathfrak{so}(2n,\C)$ has two different roots of level $k-2$ which are
$\rho=\alpha_1+\alpha_2+\alpha_3+2\sum_{i=4}^{n}\alpha_i+\alpha_{n-1}+\alpha_n$ and $\rho'=\alpha_2+\sum_{i=3}^{n}2\alpha_i+\alpha_{n-1}+\alpha_n$.

There are two roots of level three if $n\geq6$ and in there is only one in other case; this is why we add the hypothesis $n\geq 6$ to repeat the proof made for the family $\mathfrak{so}(2n+1,\C)$. In this case the roots of level $k-3$ are $\theta=\alpha_1+\alpha_2+\alpha_3+\alpha_4+2\sum_{i=5}^{n-2}\alpha_i+\alpha_{n-1}+\alpha_n$ and $\theta'=\alpha_2+\alpha_3+2\sum_{i=4}^{n-2}\alpha_i+\alpha_{n-1}+\alpha_n$.

For these roots the same relations as for the last case hold:
$$
\alpha_{\max}=\delta +\alpha_2,\qquad
\alpha_{\max}=\rho+(\alpha_2+\alpha_3)\;=\;\rho'+(\alpha_1+\alpha_2)$$$$
\alpha_{\max}=\theta +(\alpha_2+\alpha_3+\alpha_4)\;
=\;\theta'+(\alpha_1+\alpha_2+\alpha_3).$$

Then Eq. \prettyref{eq:eq85} is valid, and so are Eqns. \prettyref{eq:eq86} and \prettyref{eq:eq87}. The roots of level two, three and four are
\begin{eqnarray}
\ell =2:&&\quad\alpha_i+\alpha_{i+1}, \,i=1,\ldots,n-2 \mbox{ and }\alpha_{n-2}+\alpha_n.\nonumber\\
\ell =3:&&\quad\alpha_i+\alpha_{i+1}+\alpha_{i+2}, \,i=1,\ldots,n-2 \mbox{ and }\alpha_{n-3}+\alpha_{n-2}+\alpha_n.\nonumber\\
\ell =4:&&\quad\alpha_i+\alpha_{i+1}+\alpha_{i+2}+\alpha_{i+3}, \,i=1,\ldots,n-3 \mbox{ and }\alpha_{n-4}+\alpha_{n-3}+\alpha_{n-2}+\alpha_n.\nonumber\end{eqnarray}
Notice that the differentials of the elements of the basis of $L_2^*,\,L_3^*$ and $L_4^*$ do not coincide with those in the previous case. Nevertheless, they have the same behavior. Proceeding in an analogous manner we obtain also in this case that $E_\infty^{0,2}(\mn)=0$.

For $n=4$ and $n=5$ we used a computational program to verify that $E_\infty^{0,2}(\mn)=0$ in both cases.
\end{Proof}

\medskip

From the proof of the classification Theorem we can state the following:
\begin{Corollary} \label{cor:cor2} For a nilpotent Lie algebra  $\mn$ as in Theorem \ref{teo:teorema} of dimension $\geq 2$ the followings conditions are equivalent:
\begin{enumerate}
\item any even dimensional trivial extension $\R^s\oplus \mn$ is symplectic,
\item $E_\infty^{0,2}( \mn)\neq 0$,
\item $\mgg=\mathfrak{sl}\,(3,\C)$ or $\mgg= \mathfrak{so}\,(5,\C).$
\end{enumerate}
\end{Corollary}

\medskip

\noindent{\em Acknowledgments.} This paper is part of my Ph.D. thesis, written at FCEIA, Universidad Nacional de Rosario, Argentina. I am grateful to my advisor Isabel Dotti for her commitment on the guidance work. I also wish to thank Jorge Lauret and Roberto Miatello for their useful suggestions.


\bibliographystyle{plain}

\bibliography{biblio}

\begin{thebibliography}{10}

\bibitem{BE-GO}
C.~Benson and C.~Gordon.
\newblock K\"ahler and symplectic structures on nilmanifolds.
\newblock {\em Topology}, 27(4):513--518, 1988.

\bibitem{dB2}
V.~del Barco.
\newblock Canonical decomposition of the cohomology groups of nilpotent {L}ie
  algebras.
\newblock arXiv:1204.4123v2, 2011.

\bibitem{dB}
V.~del Barco.
\newblock Symplectic structures on free nilpotent {L}ie algebras.
\newblock arXiv:1111.3280v1, 2011.

\bibitem{DI2}
J.~Dixmier.
\newblock Cohomologie des alg\`ebres de {L}ie nilpotentes.
\newblock {\em Acta Sci. Math. Szeged}, 16:246--250, 1955.

\bibitem{DO-TI}
I.~Dotti and P.~Tirao.
\newblock Symplectic structures on {H}eisenberg-type nilmanifolds.
\newblock {\em Manuscripta Math.}, 102(3):383--401, 2000.

\bibitem{BO-GO}
M.~Goze and A.~Bouyakoub.
\newblock Sur les alg\`ebres de {L}ie munies d'une forme symplectique.
\newblock {\em Rend. Sem. Fac. Sci. Univ. Cagliari}, 57(1):85--97, 1987.

\bibitem{HA}
S.~Halperin.
\newblock Le complexe de {K}oszul en alg\`ebre et topologie.
\newblock {\em Ann. Inst. Fourier (Grenoble)}, 37(4):77--97, 1987.

\bibitem{HE}
S.~Helgason.
\newblock {\em Differential geometry, {L}ie groups, and symmetric spaces},
  volume~80 of {\em Pure and Applied Mathematics}.
\newblock Academic Press Inc. [Harcourt Brace Jovanovich Publishers], New York,
  1978.

\bibitem{KN2}
A.~Knapp.
\newblock {\em {L}ie Groups, {L}ie algebras, and cohomology}.
\newblock Lecture notes in mathematics. Princeton University Press, 1988.

\bibitem{KO}
B~Kostant.
\newblock Lie algebra cohomology and the generalized {B}orel-{W}eil theorem.
\newblock {\em Ann. of Math. (2)}, 74:329--387, 1961.

\bibitem{MI}
D.~Millionschikov.
\newblock Graded filiform {L}ie algebras and symplectic nilmanifolds.
\newblock In {\em Geometry, topology, and mathematical physics}, volume 212 of
  {\em Amer. Math. Soc. Transl. Ser. 2}, pages 259--279. Amer. Math. Soc.,
  Providence, RI, 2004.

\bibitem{NO}
K.~Nomizu.
\newblock On the cohomology of compact homogeneous spaces of nilpotent {L}ie
  groups.
\newblock {\em Ann. of Math. (2)}, 59:531--538, 1954.

\bibitem{PO-TI}
H.~Pouseele and P.~Tirao.
\newblock Compact symplectic nilmanifolds associated with graphs.
\newblock {\em J. Pure Appl. Algebra}, 213(9):1788--1794, 2009.

\bibitem{SA1}
S.M. Salamon.
\newblock Complex structures on nilpotent {L}ie algebras.
\newblock {\em J. Pure Appl. Algebra}, 157(2-3):311--333, 2001.

\bibitem{SI}
J.~{\v{S}}ilhan.
\newblock A real analog of {K}ostant's version of the {B}ott-{B}orel-{W}eil
  theorem.
\newblock {\em J. Lie Theory}, 14(2):481--499, 2004.

\end{thebibliography}

\end{document}